\documentclass[12pt]{article}
\usepackage{color,latexsym,amsfonts,amssymb}
\usepackage{amsmath}
\usepackage{amsfonts,amssymb}
\def\nexto{\kern -0.54em}
\newcommand{\mean}{{\rm {I\ \nexto E}}}
\def\prob{{\rm {I\ \nexto P}}}
\def\var{{\rm Var}}
\def\real{{\rm {I\ \nexto R}}}

\newcommand{\bn}{{\bf N}}

\newcommand{\bz}{{\bf Z}}

\newcommand{\bl}{\mbox{\boldmath$\lambda$}}

\newcommand{\boldeta}{\mbox{\boldmath$\eta$}}

\newcommand{\bm}{\mbox{\boldmath$\mu$}}

\newcommand{\boldXi}{\mbox{\boldmath$\Xi$}}
\newcommand{\boldxi}{\mbox{\boldmath$\xi$}}

\def\bone{{\bf 1}}

\newcommand{\cb}{{\cal B}}
\newcommand{\cf}{{\cal F}}

\newcommand{\ch}{{\cal H}}

\newcommand{\cs}{{\cal S}}

\def\equald{\stackrel{\mbox{\scriptsize{{\rm d}}}}{=}}

\def\ignore#1{}

\newcommand{\qed}{\hbox{\vrule width 5pt height 5pt depth 0pt}}
\def\Ref#1{(\ref{#1})}

\def\blct{{\bl_{C,t}}}
\def\Xict{\Xi_{C,t}}
\def\bXict{\boldXi_{C,t}}

\def\tXict{{\tilde \Xi_{C,t}}}
\def\Z{\mathbb{Z}}

\newcommand{\beas}{\begin{eqnarray*}}
\newcommand{\enas}{\end{eqnarray*}}

\newcommand{\bea}{\begin{eqnarray}}
\newcommand{\ena}{\end{eqnarray}}

\def\eq{\begin{equation}}
\def\en{\end{equation}}

\newtheorem{Theorem}{Theorem}[section]

\newtheorem{Lemma}[Theorem]{Lemma}

\newtheorem{Proposition}[Theorem]{Proposition}
\newtheorem{Remark}[Theorem]{Remark}

\makeatletter
\renewcommand\theequation{\thesection.\@arabic\c@equation}
\makeatother

\begin{document}

\title{On customer flows in Jackson queuing networks}
\author{Sen Tan\ and\
Aihua Xia\footnote{Corresponding author. 
E-mail: xia@ms.unimelb.edu.au}\\
Department of Mathematics and Statistics\\
The University of
Melbourne\\
Parkville, VIC 3052\\
Australia}

\date {{19 January, 2010}}

\maketitle

\begin{abstract} Melamed's theorem states that for a Jackson queuing network, the equilibrium flow along a link 
follows Poisson distribution if and only if no customers
can travel along the link more than once. Barbour \& Brown~(1996) considered the Poisson approximate version of Melamed's theorem by
allowing the customers a small probability $p$ of travelling along the link more than once. In this paper, we prove that the customer flow process is a Poisson cluster process and then establish a general approximate version of Melamed's theorem accommodating all possible cases of  
$0\le p<1$.


\end{abstract}

\section{Introduction}
\setcounter{equation}{0}

We consider a Jackson queuing network 
with $J$ queues and the following specifications [see Barbour \& Brown~(1996) for more details]. First, we assume that customers can move from one queue to another as well as can enter and leave from any queue. We assume that the exogenous arrival processes are independent Poisson processes with rates $\nu_j$, $1\le j\le J$. Service requirements are assumed to be exponential random variables with parameter 1 and when there are $m$ customers in queue $j$,
the service effort for queue $j$ is $\phi_j(m)$, where $\phi_j(0)=0$, $\phi_j(1)>0$ and  $\phi_j(m)$ is a non-decreasing function of $m$. 
Second, we define the switching process as follows. 
Let $\lambda_{ij}$ be the probability that an individual moves from queue $i$ to queue $j$,
 $\mu_i$ be the exit probability from queue $i$ and it is natural to assume
 $$\sum_{j=1}^J\lambda_{ij}+\mu_i=1, \ 1\le i\le J.$$
 Without loss of generality, we may assume that the network is irreducible in the sense that all customers can access any queue with a positive probability.
 Set $\alpha_j$  as the total rate of arriving customers (including both exogenous and endogenous arrivals) to queue $j$, then the rates $\{\alpha_j\}$ satisfy 
the equations 
$$\alpha_j=\nu_j+\sum_{i=1}^J\alpha_i\lambda_{ij},\ 1\le j\le J$$
and they are the unique solution of the equations with $\alpha_j>0$ for all $j$.

For convenience, we define state 0 as the outside of the network, that is, the point of arrival and departure of an individual into and from the system.
We write $\cs:=\{(j,k):\ 0\le j,k\le J\}$ as the set of all possible direct links and use
 $\boldXi^{jk}$ to record the transitions of individuals moving from queue $j$ to queue $k$, then
 $\boldXi=\{\boldXi^{jk},\ 0\le j,k\le J\}$ gives a full account of customer flows in the network, where departures are transitions to 0 and arrivals are transitions from 0. If $\rho_{jk}$ is the rate of equilibrium flow along the link $(j,k)$, then  $\rho_{jk}=\alpha_j\lambda_{jk}$ and the mean measure of $\boldXi$
 is
 $$\bl(ds,(j,k))=\rho_{jk}ds,\ s\in\real,\ (j,k)\in \cs.$$

 Our interest is on the customer flows along the links in $C\subset{\cal S}$ for the time interval $[0,t]$, so we set the carrier space as $\Gamma_{C,t}=[0,t]\times C$ and use $\bXict$ to stand for the transitions along the links in $C$ for the period $[0,t]$. Then the mean measure of $\bXict$
 is
 $$\blct(ds,(j,k))=\rho_{jk}ds,\ 0\le s\le t,\ (j,k)\in C.$$

Melamed's theorem states that $\bXict$ is a Poisson process if and only if no customers travel along the links in $C$ more than once [Melamed~(1979) proved the theorem when $\phi_j(m)=c_j$ for $m\ge 1$ and $1\le j\le J$, and the general case was completed by Walrand \& Varaiya~(1981)]. Barbour \& Brown~(1996) considered the Poisson approximate version of Melamed's theorem by
allowing the customers a small probability of traveling along the links more than once. For convenience, we call the probability of customers traveling along the links in $C$ more than once as the {\it loop probability}. The bounds for the errors of Poisson process approximation are sharpened by Brown, Weinberg \& Xia~(2000) and Brown, Fackrell \& Xia~(2005) and it is concluded in these studies that the accuracy of Poisson approximation depends on how small the loop probability is. 

In section~2 of the paper, we use the Palm theory, Barbour-Brown Lemma [Barbour \& Brown~(1996)] and infinite divisibility of point processes to prove that 
$\bXict$ is a Poisson cluster process [Daley \& Vere-Jones~(1988), p.~243]. The characterization involves a few quantities which are generally intractable, so in section~3, we prove that $\boldXi$ is over-dispersed [see Brown, Hamza \& Xia~(1998)], i.e., its variance is greater than its mean, and conclude that suitable approximations should be those with the same property, such as the compound Poisson or negative binomial distributions. We then establish a general approximate version of Melamed's theorem for the total number of customers traveling along the links in $C$ for the period $[0,t]$ based on a suitably chosen negative binomial distribution. The approximation error is measured in terms of the total variation distance and the error bound is small when the loop probability is small [cf. the Poisson approximation error bound in Barbour \& Brown~(1996)] and/or $t$ is large with order $\frac{1}{\sqrt{t}}$ [cf. Berry--Esseen bound for normal approximation, Petrov~(1995)].

\section{A characterization of the customer flow process}
\setcounter{equation}{0}

The customer flows are directly linked to the changes of the states of the queue lengths, so we define $N_i(t)$ as the number of customers, including those in service, at queue $i$ and time $t$, $1\le i\le J$. Then
 $\{\bn(t):=(N_1(t),\dots,N_J(t)):\ t\in\real\}$ is a pure Markov jump process on state space $\{0,1,2,\dots\}^J$ with the following transition rates ${\bf n}=(n_1,\dots,n_J)$:
$$q_{{\bf n},{\bf n}+e_j}=\nu_j,\ q_{{\bf n},{\bf n}-e_j+e_k}=\phi_j(n_j)\lambda_{jk},\ q_{{\bf n},{\bf n}-e_j}=\mu_j\phi_j(n_j),\ 1\le j,k\le J,$$
where $e_j$ is the $j$th coordinate vector in $\{0,1,2,\dots\}^J$.
The stationary Markov queue length process has a unique stationary distribution: for each $t>0$, $N_j(t)$, $j=1,\dots,J$ are independent with 
$$\prob(N_j(t)=k)=\frac{\alpha_j^k/\prod_{r=1}^k\phi_j(r)}{\sum_{l=0}^\infty
\alpha_j^l/\prod_{r=1}^l\phi_j(r)}.$$

Let $X_i$ be the $i$th queue visited by a given customer, then
 $\{X_i:\ i=1,2,\dots\}$ is called the {\it forward customer chain} and it is a homogeneous finite Markov chain with transition probabilities
\beas
&&p_{00}=0,\ p_{0k}=\frac{\nu_k}{\sum_{j=1}^J\nu_j};\\
&&p_{j0}=\mu_j,\ p_{jk}=\lambda_{jk},\ j,k=1,\dots,J.
\enas
The backward customer chain $X^\ast$ is the forward customer chain for the time-reversed process of $\{\bn(t):=(N_1(t),\dots,N_J(t)):\ t\in\real\}$ [Barbour \& Brown~(1996), p.~475] and it can be viewed as the time-reversal of the forward customer chain $\{X_i:\ i=1,2,\dots\}$ with transition probabilities
\beas
&&p_{00}^\ast=0,\ p_{0j}^\ast=\frac{\mu_j\alpha_j}{\sum_{l=1}^J\mu_l\alpha_l};\\
&&p_{k0}^\ast=\frac{\nu_k}{\alpha_k},\ p_{kj}^\ast=\frac{\alpha_j\lambda_{jk}}{\alpha_k},\ j,k=1,\dots,J.
\enas

We will use the Palm distributions to characterize the distribution of $\bXict$, prove its properties and establish 
a general approximate version of Melamed's theorem for $\bXict(\Gamma_{C,t})$. 
For the point process $\boldXi$ with locally finite mean measure $\bl(ds,(j,k))=\rho_{jk}ds,\ s\in\real,\ (j,k)\in\cs$, we may consider it as a random measure on the metric space $\real\times \cs$ equipped with the metric 
$$d((u_1,(j_1,j_2)),(u_2,(k_1,k_2)))=|u_1-u_2|\bone_{(j_1,j_2)\ne(k_1,k_2)},\ u_1,u_2\in\real\mbox{ and }(j_1,j_2),\ (k_1,k_2)\in\cs,$$ so that we can define the {\it Palm distribution} at $\alpha\in \real\times \cs$ as the distribution of $\boldXi$ conditional on the presence of a point at $\alpha$, that is,
$$P^\alpha(\cdot)=
\frac{\mean \left[1_{[{\scriptsize \boldXi}\in\cdot]}\boldXi(d\alpha)\right]}{\bl(d\alpha)},\
\alpha\in\Gamma_{C,t}\
\bl-\mbox{almost surely},$$
see Kallenberg~(1983), p.~83 for more details. 
A process $\boldXi^\alpha$ is called the {\it Palm process} of $\boldXi$ at $\alpha$ if its distribution is $P^\alpha$. In applications, it is often more convenient to work with the {\it reduced Palm process} $\boldXi^\alpha-\delta_\alpha$
[Kallenberg~(1983), p.~84], where $\delta_\alpha$ is the Dirac measure at $\alpha$.
 
The Palm distributions are closely related to the {\it size-biasing} in 
sampling contexts [Cochran~(1977)]. More precisely, if $X$ is a non-negative integer-valued random variable, one may consider it as a point process with the carrier space having only one point so its Palm distribution becomes 
$$\prob(X_s=i):=\frac{i\prob(X=i)}{\mean X}.$$ 
However, this is exactly the definition of the {\it size
biased} distribution of $X$ [see Goldstein \& Xia~(2006)].

\begin{Lemma}\label{BBlemma}
[Barbour \& Brown (1996)]
For the open queuing network, the reduced Palm distribution for the network given a transition at link 
$(j,k)$ at time $0$ is the same as that for the original network, save that the network on $(0,\infty)$ behaves as if there were an extra individual at queue $k$ at time 0 and the network on $(-\infty,0)$ behaves as if there were an extra individual in queue $j$ at time 0.
\end{Lemma}

For two random elements $\eta_1$ and $\eta_2$ having the same distribution, we write for brevity $\eta_1\equald\eta_2$.

\begin{Lemma}\label{XiaLemma1} For each $(j,k)\in\cs$, there is a point process $\boldxi^{(0,(j,k))}$ on $\real\times\cs$ independent of $\boldXi$ such that 
$$\boldxi^{(0,(j,k))}+\boldXi \equald\boldXi^{(0,(j,k))}.$$
\end{Lemma}

\noindent{\bf Proof.} The proof is adapted from Barbour \& Brown~(1996), p.~480. By Lemma~\ref{BBlemma}, the reduced Palm process $\boldXi^{(0,(j,k))}-\delta_{(0,(j,k))}$ has the same distribution as that of $\boldXi$ except that the network on $(0,\infty)$ behaves as if there were an extra individual at queue $k$ at time 0 and the network on $(-\infty,0)$ behaves as if there were an extra individual in queue $j$ at time 0. Let $\tilde X^{(0)}$ and $X^{(0)}$ be the routes taken by the extra individual on $(-\infty,0)$ and $(0,\infty)$ respectively. Whenever the extra customer is at queue $i$ together with other $m$ customers, we use independently sampled exponential service requirements with instantaneous service rate $\phi_i(m+1)-\phi_i(m)$. Noting that this construction ensures that the extra customer uses the ``spare" service effort and never ``interferes" with the flow of the main traffic, one can see that its transitions are independent of $\boldXi$. The same procedure applies to the construction of the backward route. Let
$\boldxi^{(0,(j,k))}$ be the transitions taken by the extra customer on $(-\infty,0)\cup(0,\infty)$ plus the Dirac measure $\delta_{(0,(j,k))}$, then $\boldxi^{(0,(j,k))}$ is independent of $\boldXi$ and the conclusion of the lemma follows from the construction. \qed

Let $\theta_s, \ s\in\real$, denote the shift operator on  
$\real\times\cs$ which translates each point in $\real\times\cs$ by $s$ to the left, i.e. $\theta_s((u,(j,k)))=(u-s,(j,k))$ and we use 
$\boldxi^{(s,(j,k))}$ to stand for a copy of $\boldxi^{(0,(j,k))}\circ\theta_s$, $s\in\real$.

From now on, we focus on the point process $\bXict$. With metric $d$, $\Gamma_{C,t}$ is a Polish space and we use $\cb\left(\Gamma_{C,t}\right)$ to stand for the Borel $\sigma$-algebra in $\Gamma_{C,t}$. Let $H_{C,t}$ denote the class of all configurations  (finite nonnegative integer-valued measures) on $\Gamma_{C,t}$ 
with $\ch_{C,t}$ the $\sigma$-algebra in $H_{C,t}$ generated by the sets 
$$\{\xi\in H_{C,t}: \xi(B)=i\},\ i\in\Z_+:=\{0,1,2,\dots\},
\ B\in\cb\left(\Gamma_{C,t}\right),$$
see Kallenberg~(1983), p.~12.

\begin{Theorem}\label{Th1} Let $\{\boldeta_i,\ i\ge 0\}$ be independent and identically distributed random measures on $\Gamma_{C,t}$ having the distribution
\eq\prob\left[\boldeta_0\left(\Gamma_{C,t}\right)\ge 1\right]=1,\ \prob(\boldeta_0\in A)=\mean \sum_{(j,k)\in C}\int_0^t
\frac{\bone_{[{\scriptsize\boldxi}^{(s,(j,k))}\in A]}}{\boldxi^{(s,(j,k))}\left(\Gamma_{C,t}\right)}\cdot\frac{\rho_{jk}}{\theta_{C,t}}ds,\ A\in \ch_{C,t},
\label{Xia2.1}
\en
where 
\eq\theta_{C,t}=\mean \sum_{(j,k)\in C}\int_0^t
\frac{1}{\boldxi^{(s,(j,k))}\left(\Gamma_{C,t}\right)}\rho_{jk}ds.\label{Xia2.2}\en
Let $M$ be a Poisson random variable with mean $\theta_{C,t}$ and independent of $\{\boldeta_i,\ i\ge 0\}$, then
$$\bXict\equald \sum_{i=1}^M\boldeta_i.$$
\end{Theorem}

\noindent{\bf Proof.} By Lemma~\ref{XiaLemma1} and Theorem~11.2 of [Kallenberg~(1983)], we can conclude that $\bXict$ is infinitely divisible, hence we obtain from Lemma~6.6 and Theorem~6.1 of [Kallenberg~(1983)] that $\bXict$ is a Poisson cluster process, that is, 
$$\bXict\equald \sum_{i=1}^M\boldeta_i,$$
where $\boldeta_i,\ i\ge 0$ are independent and identically distributed random measures on $\Gamma_{C,t}$ such that $\prob\left(\boldeta_0(\Gamma_{C,t})\ge 1\right)=1$, $M$ is a Poisson random variable with mean $\theta_{C,t}$ and independent of  $\{\boldeta_i,\ i\ge 1\}$.
The direct verification ensures that the Palm process of $\sum_{i=1}^M\boldeta_i$ at $\alpha\in\Gamma_{C,t}$ is 
$\sum_{i=1}^M\boldeta_i+\boldeta_0^{\alpha}$, where $\boldeta_0^{\alpha}$ is the Palm process of $\boldeta_0$ at $\alpha$ and is independent of $\{M,\boldeta_i,\ i\ge 1\}$. This in turn implies that $\boldxi^{(s,(j,k))}\equald \boldeta_0^{(s,(j,k))}.$ 

Let $\bm(ds,(j,k))$ denote the mean measure of the point process $\boldeta_0$, then some elementary computation ensures that the mean measure of $\sum_{i=1}^M\boldeta_i$ is $\theta_{C,t}\bm(ds,(j,k))$ for $(j,k)\in C$ and $0\le s\le t$. On the other hand, the mean measure of $\bXict$ is $\blct(ds,(j,k))=\rho_{jk}ds $, $(j,k)\in C$ and $s\in[0,t]$, so we obtain 
\eq\bm(ds,(j,k))=\frac{\rho_{jk}}{\theta_{C,t}}ds,\ (j,k)\in C,\ s\in[0,t].\label{Xia2.3}\en
The representation \Ref{Xia2.1} is because of the fact that
$\prob\left(\boldeta_0(\Gamma_{C,t})\ge 1\right)=1$ and 
$$\prob\left(\boldeta_0\in A\right)=\mean\int_{\Gamma_{C,t}}\frac{\bone_{[{\scriptsize\boldeta_0}\in A]}}{\boldeta_0\left(\Gamma_{C,t}\right)}\boldeta_0(d\alpha)=\mean \sum_{(j,k)\in C}\int_0^t
\frac{\bone_{[{\scriptsize\boldxi}^{(s,(j,k))}\in A]}}{\boldxi^{(s,(j,k))}\left(\Gamma_{C,t}\right)}\frac{\rho_{jk}}{\theta_{C,t}}ds.$$
In particular, if we take $A=\ch_{C,t}$, then the left hand side becomes 1, so \Ref{Xia2.2} follows.
\qed

Despite the fact that $\theta_{C,t}$ is specified by \Ref{Xia2.2}, since the Palm process $\boldxi^{(s,(j,k))}$ is generally intractable, it is virtually impossible to express $\theta_{C,t}$ explicitly in terms of the specifications of the Jackson queuing network. On the other hand, the relationship \Ref{Xia2.3} yields
$$\mean\boldeta_0(\Gamma_{C,t})=\frac{\rho_C}{\theta_{C,t}}t.$$
The following proposition tells us the range of values that  $\mean\boldeta_0(\Gamma_{C,t})$ and $\theta_{C,t}$ may take. To this end, we define
\eq\epsilon_C(j,k)=\mean\boldxi^{(0,(j,k))}(\real\times C)-1\mbox{ and }\epsilon_C=\sum_{(j,k)\in C}\frac{\rho_{jk}}{\rho_C}\epsilon_C(j,k).\label{Xia2.6}\en
In other words, $\epsilon_C(j,k)$ is the average number of visits in $C$ by the extra customer crossing the link $(j,k)$ and $\epsilon_C$ is the weighted average number of visits by an extra customer crossing links in $C$. 

\begin{Proposition} We have
\eq1\le \mean\boldeta_0(\Gamma_{C,t})\le 1+\epsilon_C\label{Xia2.4}\en
and
\eq\frac{\rho_C}{1+\epsilon_C}t\le\theta_{C,t}\le\rho_C t.\label{Xia2.5}\en
\end{Proposition}

\noindent{\bf Proof.} The first inequality of \Ref{Xia2.4} follows immediately from the fact that \\
$\prob(\boldeta_0(\Gamma_{C,t})\ge 1)=1$. For the second inequality of \Ref{Xia2.4}, noting that the mean measure of $\boldeta_0$ is
$$\bm(ds,(j,k))=\frac{\rho_{jk}}{\theta_{C,t}}ds,\ (j,k)\in C,\ s\in[0,t],$$
we have
\beas
\left\{\frac{\rho_C}{\theta_{C,t}}t\right\}^2&=&[\mean\boldeta_0(\Gamma_{C,t})]^2\le\mean\left[\boldeta_0(\Gamma_{C,t})^2\right]\\
&=&\mean\int_{\Gamma_{C,t}}\boldeta_0(\Gamma_{C,t})\boldeta_0(d\alpha)\\
&=&\sum_{(j,k)\in C}\int_0^t\mean \boldeta_0^{(s,(j,k))}(\Gamma_{C,t})\bm(ds,(j,k))\\
&\le&\sum_{(j,k)\in C}(1+\epsilon_C(j,k))\frac{\rho_{jk}}{\theta_{C,t}}t\\
&=&\frac{\rho_C}{\theta_{C,t}}t+\frac{\sum_{(j,k)\in C}\epsilon_C(j,k)\rho_{jk}}{\theta_{C,t}}t\\
&=&\frac{\rho_C}{\theta_{C,t}}t(1+\epsilon_C).
\enas
We divide both sides by $\frac{\rho_C}{\theta_{C,t}}t$ to get 
$$\mean\boldeta_0(\Gamma_{C,t})=\frac{\rho_C}{\theta_{C,t}}t\le 1+\epsilon_C.$$
Finally, \Ref{Xia2.5} is an immediate consequence of \Ref{Xia2.4} and the equation $\theta_{C,t}=\frac{\rho_C}{\mean\boldeta_0(\Gamma_{C,t})}t$. \qed

\begin{Remark}{\rm If the loop probability in $C$ is 0, then $\epsilon_C=0$ and $\mean\boldeta_0(\Gamma_{C,t})=1$, so there is only one customer on $\Gamma_{C,t}$. This customer is crossing the link $(j,k)$ with probability $\frac{\rho_{jk}}{\rho_C}$ at a time uniformly distributed on $[0,t].$}
\end{Remark}

\section{A discrete central limit theorem for the customer flow process}
\setcounter{equation}{0}

A random variable is said to be {\it over-dispersed} (resp. {\it under-dispersed}) if its variance to mean ratio is greater (resp. less) than one.
A random measure $\chi$ on a Polish space is said to be {\it over-dispersed} (resp. {\it under-dispersed}) if $\chi(B)$ is over-dispersed (resp. under-dispersed) for all bounded Borel subset $B$ of the Polish space. 
It is concluded in Brown, Hamza \& Xia~(1998) that point processes which arise from Markov chains which are
time-reversible, have finitely many states and are irreducible are always
over-dispersed. As our process $\boldXi$ is virtually a multivariate version of point processes studied in Brown, Hamza \& Xia~(1998), the following property can be viewed as a natural extension of the study in Brown, Hamza \& Xia~(1998).


\begin{Proposition} The point process $\boldXi$ is over-dispersed.
\end{Proposition}

\noindent{\bf Proof.} The space $(\real\times\cs,d)$ is a Polish space and for each bounded Borel subset $B$ of $\real\times \cs$, it follows from the definition of the Palm processes [see Kallenberg~(1983), p.~84, equation~(10.4)] that
\beas
&&\mean\left[\boldXi(B)\right]^2=\mean\int_B\boldXi(B)\boldXi(d\alpha)=\mean \int_B\boldXi^\alpha(B)\bl(d\alpha)
\ge\mean \int_B\left(\boldXi(B)+1\right)\bl(d\alpha),
\enas
that is,
\eq\var\left[\boldXi(B)\right]\ge\mean \boldXi(B),\label{xia3.1}
\en
completing the proof. \qed

The inequality in \Ref{xia3.1} is generally strict except that the loop probability is 0, 
i.e. $\boldXi$ is a Poisson process.
Hence, suitable approximate models for the distribution of $\Xict:=\bXict(\Gamma_{C,t})$
are necessarily over-dispersed. One potential candidate for approximating the distribution of $\Xict$ is the compound Poisson distribution. However, as it is virtually impossible to extract the distribution of
$\boldxi^{(s,(j,k))}(\Gamma_{C,t})$ for $0\le s\le t$, we face the same difficulty to specify and estimate the approximate distribution if a general compound Poisson is used.  On the other hand, as a special family of the compound Poisson distributions [Johnson, Kemp \& Kotz (2005), pp.~212--213 and p.~346], the negative binomial distribution has been well documented as a natural model for many over-dispersed random phenomena [see Bliss \& Fisher~(1953), Wang \& Xia~(2008) and Xia \& Zhang~(2009)]. The negative binomial distribution
${\rm NB}(r,q)$, $r>0$, $0<q<1$, is defined as
$$\pi_i=\frac{\Gamma(r+i)}{\Gamma(r)i!}q^r(1-q)^i,\ i\in\Z_+.$$
The advantage of using negative binomial approximation is that it suffices to estimate the mean and variance of the approximating distribution, like what we often do in applying the central limit theorem based on the normal approximation.

We will use the total variation distance between the distributions of nonnegative integer-valued random variables $Y_1$ and $Y_2$ 
$$d_{TV}(Y_1,Y_2):=\sup_{A\subset  \Z_+}|\prob(Y_1\in A)-\prob(Y_2\in A)|$$
to measure the approximation errors in negative binomial approximation.

The discrete central limit theorem is valid under the assumption that the loop probability  
is less than 1. More precisely, let $w_C(jk)$ be the probability that a link $(j,k)$ crossing customer crosses the links in $C$ only once, i.e., the only time that the customer crosses the links in $C$ is the one the customer is crossing. Define
$$w_C=\sum_{(j,k)\in C}w_C(jk)\rho_{jk}/\rho_C,$$
the weighted probability of customers crossing links in $C$ only once. Clearly, we have $w_C(jk)\ge\mu_k$, so
$$w_C\ge\sum_{(j,k)\in C}\rho_{jk}\mu_k/\rho_C.$$ 
The following lemma plays a crucial rule for the estimation of the negative binomial approximation error.

\begin{Lemma}\label{keylemma1}
$d_{TV}\left(\Xict,\Xict+1\right)\le \frac1{\sqrt{2e w_C\rho_Ct}}.$
\end{Lemma}

\noindent{\bf Proof.} We prove the claim by a coupling based on the ``priority principle" [cf. the proof of Lemma~\ref{XiaLemma1}]. We define a customer as a {\it single crossing} ({\it sc} for brevity) customer if the customer crosses links in $C$ only once, otherwise, the customer is labeled as {\it multiple crossing}, or {\it mc} for short. We ``manage" the network by regrouping the customers at each queue into {\it sc} customers and {\it mc} customers. Whenever there are $m_2$ {\it mc} customers together with $m_1$ {\it sc} customers at queue $j$, we use independently sampled exponential service requirements with instantaneous service rate $\phi_j(m_1+m_2)-\phi_j(m_1)$ for all of the {\it mc} customers while the service for the {\it sc} customers is carried out with instantaneous service rate $\phi_j(m_1)$, that is, as if no {\it mc} customers present at the queue. Since the {\it sc} customers take priority over the {\it mc} customers and the {\it mc} customers use the ``spare" service effort and never interrupt the traffic flow of the {\it sc} ones, we can see that its transitions are independent of the transitions of the {\it sc} customers. Let $Z_1^{jk}$ (resp. $Z_2^{jk}$) denote the transitions of {\it sc} (resp. {\it mc}) customers moving from queue $j$ to queue $k$ in the period $[0,t]$, then $\bz_1:=\{Z_1^{jk},\ (j,k)\in C\}$ and 
$\bz_2:=\{Z_2^{jk},\ (j,k)\in C\}$ are independent and
$$\bXict\equald \bz_1+\bz_2.$$
By Melamed's theorem, the point process $\bz_1$ is a Poisson process with mean measure
$$\bl_{\bz_1}(ds,(j,k))=w_C(jk)\rho_{jk}ds,\ (j,k)\in C,\ 0\le s\le t,$$
so $\bz_1(\Gamma_{C,t})$ follows Poisson distribution with mean $w_C\rho_Ct$ and
$$d_{TV}\left(\Xict,\Xict+1\right)\le d_{TV}\left(\bz_1(\Gamma_{C,t}),\bz_1(\Gamma_{C,t})+1\right)\le \frac1{\sqrt{2e w_C\rho_Ct}},$$
where the last inequality is because of the fact that the distribution of Poisson is unimodal and Proposition~A.2.7 of [Barbour, Holst \& Janson~(1992), p.~262].
\qed

To state the discrete central limit theorem, we set 
$$\sigma_C(j,k)=\mean[\boldxi^{(0,(j,k))}(\real\times C)(\boldxi^{(0,(j,k))}(\real\times C)-1)]\mbox{ and }\sigma_C=\sum_{(j,k)\in C}\frac{\rho_{jk}}{\rho_C}\sigma_C(j,k).$$
That is, $\sigma_C(j,k)$ is the second factorial moment of the number of visits in $C$ by the extra customer crossing the link $(j,k)$ and $\sigma_C$ is the weighted average of the second factorial moments of the number of visits by an extra customer crossing links in $C$ [cf. \Ref{Xia2.6}].

\begin{Theorem}\label{Xia3.2} Let 
$$r=\frac{(\rho_Ct)^2}{\var(\Xict)-\rho_Ct},\ q=\frac{\rho_Ct}{\var(\Xict)},$$
then
\bea
d_{TV}\left(\Xict,{\rm NB}(r,q)\right)&\le&\frac1{(\rho_Ct)^2\sqrt{2e w_C\rho_Ct}}
\left\{2(\var(\Xict)-\rho_Ct)^2\right.\nonumber\\
&&\mbox{\hskip0cm}\left.+\rho_Ct(\Xict[3]-\rho_Ct\Xict[2]-2\rho_Ct(\var(\Xict)-\rho_Ct))\right\}\label{Xia3.2.1}\\
&\le&\frac1{\sqrt{2e w_C\rho_Ct}}(2\epsilon_C^2+\sigma_C),\label{Xia3.2.2}
\ena
where $\Xict[n]$ stands for the $n$th factorial moment of $\Xict$ defined as
$$\Xict[n]=\mean[\Xict(\Xict-1)\dots(\Xict-n+1)].$$ 
\end{Theorem}

\begin{Remark}{\rm The parameters of the approximating negative binomial distribution are chosen so that it matches the mean and variance of $\Xict$.}
\end{Remark}

\begin{Remark}{\rm If 
the loop probability in $C$ is 0, then the negative binomial is reduced to Poisson distribution and the upper bound in Theorem~\ref{Xia3.2} becomes 0. This implies half of Melamed's theorem~(1979).}
\end{Remark}

\begin{Remark}{\rm
If the loop probability 
is between 0 and 1, then both $\epsilon_C$ and $\sigma_C$ are finite, so the negative binomial approximation error bound is of order $O(1/\sqrt{t})$. Furthermore, if the loop probability is small, then both $\epsilon_C$ and $\sigma_C$ are small, so the negative binomial approximation to the distribution of $\Xict$ is even faster.}
\end{Remark}

\noindent{\bf Proof of Theorem~\ref{Xia3.2}.} The essence of Stein's method is to find a generator which characterizes the approximating distribution, establish a Stein identity to transform the problem of estimating the approximation errors into the study of the structure of the object under investigation. 
In the context of negative binomial approximation, 
let $a=r(1-q)$, $b=1-q$, then a generator which characterizes ${\rm NB}(r,q)$ is defined as
$$\cb g(i)=(a+bi)g(i+1)-ig(i),\ i\in\Z_+,$$
for all bounded functions $g$ on $\Z_+$ [see Brown \& Phillips~(1999) and Brown \& Xia~(2001)]. The Stein identity is naturally established as
\eq\cb g(i)=f(i)-\pi(f)\label{Steinidentity}\en
for $f\in\cf:=\{f:\ \Z_+\to[0,1]\}$, where $\pi(f)=\sum_{i=0}^\infty f(i)\pi_i$. It was shown in Brown \& Xia~(2001) that, for each $f\in\cf$, the solution $g_f$ to the Stein equation \Ref{Steinidentity} satisfies
\eq\|\Delta g_f\|\le\frac1a,\label{Steinconstant}\en
where $\Delta g_f(\cdot)=g_f(\cdot+1)-g_f(\cdot).$ The Stein identity \Ref{Steinidentity} ensures that
$$\sup_{f\in\cf}\left |\mean f(\Xict)-\pi(f)\right|=\sup_{f\in\cf}\left|\mean \cb g_f(\Xict)\right|,$$
hence, it suffices to estimate $\mean \cb g_f(\Xict)$ for all $f\in\cf$. For convenience, we drop $f$ from the subindex of $g_f$.
By Lemma~\ref{XiaLemma1}, we can take a point process $\boldxi_{C,t}^{(s,(j,k))}$ on $\Gamma_t$ independent of $\bXict$ such that 
$$\bXict^{(s,(j,k))}=\bXict+\boldxi_{C,t}^{(s,(j,k))}.$$
Therefore, if we write $\boldxi_{C,t}^{(s,(j,k))}(\Gamma_{C,t})=1+\xi^{(s,(j,k))}$, then
\bea
\mean\cb g(\Xict)&=&\mean[(a+b\Xict)g(\Xict+1)-\Xict g(\Xict)]\nonumber\\
&=&a\mean g(\Xict+1)+b\sum_{(j,k)\in C}\int_0^tg(\Xict+2+\xi^{(s,(j,k))})\rho_{jk}ds\nonumber\\
&&-\sum_{(j,k)\in C}\int_0^tg(\Xict+1+\xi^{(s,(j,k))})\rho_{jk}ds.\label{proof01}
\ena
Let
\eq a+(b-1)\sum_{(j,k)\in C}\rho_{jk}t=0,\label{coefficient1}\en
and $\tXict=\Xict+1$, then it follows from \Ref{proof01} that
\bea
&&\mean\cb g(\Xict)\nonumber\\
&&=\mean \sum_{(j,k)\in C}\int_0^t\left[b\left(g\left(\tXict+1+\xi^{(s,(j,k))}\right)-g\left(\tXict\right)\right)-\left(g\left(\tXict+\xi^{(s,(j,k))}\right)-g\left(\tXict\right)\right)\right]\rho_{jk}ds\nonumber\\
&&=\mean \sum_{(j,k)\in C}\int_0^t\left\{\sum_{r=0}^{\xi^{(s,(j,k))}-1}\left[b\Delta g\left(\tXict+r+1\right)
-\Delta g\left(\tXict+r\right)\right]+b\Delta g\left(\tXict\right)\right\}\rho_{jk}ds.
\label{proof02}\nonumber\\
\ena
Now, set
\eq  b=\frac{\sum_{(j,k)\in C}\int_0^t\mean \xi^{(s,(j,k))}\rho_{jk}ds}{\sum_{(j,k)\in C}\int_0^t\mean \xi^{(s,(j,k))}\rho_{jk}ds+\rho_C t}
=\frac{\var\left(\Xict\right)-\rho_Ct}{\var\left(\Xict\right)}\label{coefficient2},\en
where the last equality is due to the following observation:
\beas
\mean\Xict^2&=&\mean\int_{\Gamma_{C,t}}\Xict\bXict(d\alpha)\\
&=&\sum_{(j,k)\in C}\mean\int_0^t\left(\Xict+1+\xi^{(s,(j,k))}\right)\rho_{jk}ds\\
&=&(\mean\Xict)^2+\rho_Ct+\mean\sum_{(j,k)\in C}\int_0^t\xi^{(s,(j,k))}\rho_{jk}ds,
\enas
and so
\eq\mean\sum_{(j,k)\in C}\int_0^t\xi^{(s,(j,k))}\rho_{jk}ds=\var(\Xict)-\rho_Ct.\label{proof04}\en
We then obtain from \Ref{proof02} that 
\bea
&&\mean\cb g(\Xict)\nonumber\\
&&=\mean \sum_{(j,k)\in C}\int_0^t\left\{\sum_{r=0}^{\xi^{(s,(j,k))}-1}\left[b\Delta^2 g\left(\tXict+r\right)-(1-b)\sum_{l=0}^{r-1}\Delta^2 g\left(\tXict+l\right)\right]\right\}\rho_{jk}ds\nonumber\\
&&=\mean \sum_{(j,k)\in C}\int_0^t\left\{\sum_{r=0}^{\xi^{(s,(j,k))}-1}\left[b\mean\Delta^2 g\left(\tXict+r\right)-(1-b)\sum_{l=0}^{r-1}\mean\Delta^2 g\left(\tXict+l\right)\right]\right\}\rho_{jk}ds,
\label{proof03}\nonumber\\
\ena
where the last equation is due to the fact that $\xi^{(s,(j,k))}$ is independent of $\tXict$. On the other hand, using \Ref{Steinconstant}, we have
$$\left|\mean \Delta^2 g(\tXict+l)\right|\le 2\|\Delta g\|d_{TV}(\Xict,\Xict+1)\le \frac{2d_{TV}(\Xict,\Xict+1)}{a},$$
so it follows from \Ref{proof03} that
\bea&&\left|\mean\cb g\left(\Xict\right)\right|\nonumber\\
&&\le\frac{d_{TV}\left(\Xict,\Xict+1\right)}{a}\sum_{(j,k)\in C}\int_0^t\left[2b\mean\xi^{(s,(j,k))}+(1-b)\mean \xi^{(s,(j,k))}\left(\xi^{(s,(j,k))}-1\right)\right]\rho_{jk}ds.
\nonumber\\
\label{proof05}\ena
Using the Palm distributions of $\bXict$ together with \Ref{proof04}, we get
\beas
\Xict[3]&=&\mean\int_{\Gamma_{C,t}}(\Xict-1)(\Xict-2)\bXict(d\alpha)\\
&=&\sum_{(j,k)\in C}\int_0^t\mean\left[\left(\Xict+\xi^{(s,(j,k))}\right)\left(\Xict+\xi^{(s,(j,k))}-1\right)\right]\rho_{jk}ds\\
&=&\rho_Ct\Xict[2]+2\rho_Ct(\var(\Xict)-\rho_Ct)+\sum_{(j,k)\in C}\int_0^t\mean \xi^{(s,(j,k))}\left(\xi^{(s,(j,k))}-1\right)\rho_{jk}ds.
\enas
This in turn ensures
\eq\sum_{(j,k)\in C}\int_0^t\mean \xi^{(s,(j,k))}\left(\xi^{(s,(j,k))}-1\right)\rho_{jk}ds=\Xict[3]-\rho_Ct\Xict[2]-2\rho_Ct(\var(\Xict)-\rho_Ct).\label{proof06}\en
Consequently, combining \Ref{proof04}, \Ref{proof06} with \Ref{proof05} gives \Ref{Xia3.2.1}. 

Finally, by the definitions of $\epsilon_C$ and $\sigma_C$, we have
$$\mean\sum_{(j,k)\in C}\int_0^t\xi^{(s,(j,k))}\rho_{jk}ds\le \epsilon_C\rho_Ct$$
and
$$\sum_{(j,k)\in C}\int_0^t\mean \xi^{(s,(j,k))}\left(\xi^{(s,(j,k))}-1\right)\rho_{jk}ds\le 
\sigma_C\rho_Ct.$$
Therefore, \Ref{Xia3.2.2} follows from \Ref{Xia3.2.1}, \Ref{proof04} and \Ref{proof06}.
\qed


\def\ac{{Academic Press}~}
\def\aap{{Adv. Appl. Prob.}~}
\def\ap{{Ann. Probab.}~}
\def\anap{{Ann. Appl. Probab.}~}
\def\jap{{J. Appl. Probab.}~}
\def\jws{{John Wiley $\&$ Sons}~}
\def\ny{{New York}~}
\def\ptrf{{Probab. Theory Related Fields}~}
\def\sp{{Springer}~}
\def\spa{{Stochastic Processes Appl.}~}
\def\sv{{Springer-Verlag}~}
\def\tpa{{Theory Probab. Appl.}~}
\def\zw{{Z. Wahrsch. Verw. Gebiete}~}



\begin{thebibliography}{Dillo 83}
\typeout{References...}
\bibliography{math}
\bibliographystyle{alpha}
\bibitem{BB96} {\sc Barbour, A. D. \& Brown, T. C.~(1996)} Approximate Versions of Melamed's Theorem. \jap\textbf{33}, 
472--489.

\bibitem{BHJ} {\sc Barbour, A. D., Holst,  L. \& Janson, S.~(1992)} {\em Poisson
Approximation.\/} Oxford Univ. Press.

\bibitem{BF53} {\sc Bliss, C. \& Fisher, R. A.~(1953)} Fitting the negative binomial distribution to biological data.
Biometrics~\textbf{9}, 174--200.
\bibitem{BFX05} {\sc Brown, T. C., Fackrell, M. \& Xia, A.~(2005)} Improved results on Poisson process approximation
in Jackson networks. {\em COSMOS Journal}~\textbf{1}, 47--55.

\bibitem{BHX98} {\sc Brown, T. C., Hamza, K. \& Xia, A.~(1998)} On the Variance to Mean Ratio for Random Variables from Markov Chains and
       Point Processes.
      \jap\textbf{35}, 303--312.

\bibitem{BrP}  {\sc Brown, T. C. \& Phillips, M. J.~(1999)}  Negative Binomial
Approximation with Stein's Method. Method. Comput. Appl. Probab.~\textbf{1:4}, 407--421.

\bibitem{BWX00} {\sc Brown, T. C., Weinberg, G.~V. \& Xia, A.~(2000)} Removing Logarithms from Poisson Process Error Bounds.
   \spa\textbf{87}, 149--165.
   
\bibitem{BX01} {\sc
Brown, T. C. \& Xia, A.~(2001)} Stein's method and birth-death processes. \ap\textbf{29}, 1373--1403.

\bibitem{Cochran} {\sc Cochran, W.~(1977)} {\em Sampling techniques.\/} Wiley.

\bibitem{DVJ}
{\sc Daley, D. J. \& Vere-Jones, D.~(1988)} {\em An Introduction to the Theory of Point Processes.\/} 
Springer-Verlag, New York.

\item\label{GX06} {\sc Goldstein, L. and Xia, A.~(2006)} Zero Biasing and a Discrete Central Limit Theorem. \ap\textbf{34}, 1782--1806.

\bibitem{JKK05}
{\sc Johnson, N. L.,  Kemp, A. \& Kotz, S.~(2005)} {\em Univariate discrete distributions.\/} Third Edition. Wiley, New York.

\bibitem{Kallenberg83} {\sc Kallenberg, O.~(1983)} {\em Random Measures.\/} \ac.

\bibitem{Mel79} {\sc Melamed, B.~(1979)} Characterizations of Poisson traffic streams in Jackson queueing networks. \aap\textbf{11}, 422--438.

\bibitem{P95} {\sc Petrov, V. V.~(1995)} {\em Limit Theorems of Probability Theory: Sequences of Independent Random
Varaibles.\/} Clarendon Press, Oxford.

\bibitem{WV81} {\sc Walrand, J. \& Varaiya, P.~(1981)} Flows in queueing networks: a martingale approach. Math. Operat. Res.~\textbf{6}, 387--404.

\item\label{WX08} {\sc Wang, X. \& Xia, A.~(2008)} On negative binomial approximation to $k$-runs. \jap\textbf{45}, 456--471.

\item\label{XZ09} {\sc Xia, A. \& Zhang, M.~(2009)} On approximation of Markov binomial distributions. Bernoulli~(accepted).

\end{thebibliography}
\end{document}